\documentclass[10pt]{amsart}

\usepackage[bookmarks=true,hyperindex,pdftex,colorlinks, citecolor=blue,linkcolor=blue, urlcolor=blue
]{hyperref}
\usepackage{color,graphicx,shortvrb}
\usepackage[dvipsnames]{xcolor}

\usepackage[active]{srcltx} 

\usepackage{enumerate}
\usepackage{amssymb}
\usepackage{tikz-cd}

\usepackage [ all ]{xy}

\newtheorem{theorem}{Theorem}[section]

\newtheorem{corollary}[theorem]{Corollary}

\newtheorem{lemma}[theorem]{Lemma}

\newtheorem{problem}[theorem]{Problem}
\newtheorem{proposition}[theorem]{Proposition}
\newtheorem{remark}[theorem]{Remark}
\newtheorem{remarks}[theorem]{Remarks}

\def\J#1#2#3{ \left\{ #1,#2,#3 \right\} }

\def\11{\textbf{$1$}}
\def\11b#1{\mathbf{1}_{_{#1}}}
\def\CC{{\mathbb{C}}}

\DeclareMathOperator{\re}{Re}

\DeclareGraphicsExtensions{.jpg,.pdf,.png,.eps}

\DeclareMathOperator{\id}{Id}

\begin{document}
	
	\title[The polynomial Daugavet property]{The Daugavet equation for polynomials on C$^*$-algebras and JB$^*$-triples}
	
	\author[D.\ Cabezas]{David Cabezas}
	\address[Cabezas]{Departamento de An\'{a}lisis Matem\'{a}tico, Facultad de
		Ciencias, Universidad de Granada, 18071 Granada, Spain.}
	\email{dcabezas@ugr.es}

	\author[M.\ Mart{\'i}n]{Miguel Mart{\'i}n}
	\address[Mart\'{\i}n]{Instituto de Matem\'{a}ticas de la Universidad de Granada (IMAG). Departamento de An\'{a}lisis Matem\'{a}tico, Facultad de	Ciencias, Universidad de Granada, 18071 Granada, Spain.}
	\email{mmartins@ugr.es}

	\author[A.\ M.\ Peralta]{Antonio M.\ Peralta}
	
	\address[Peralta]{Instituto de Matem\'{a}ticas de la Universidad de Granada (IMAG). Departamento de An\'{a}lisis Matem\'{a}tico, Facultad de	Ciencias, Universidad de Granada, 18071 Granada, Spain.}
	\email{aperalta@ugr.es}

	
	
	
	\date{June 23rd, 2022}
	
	\begin{abstract} We prove that every JB$^*$-triple $E$ (in particular, every $C^*$-algebra) satisfying the Daugavet property also satisfies the stronger polynomial Daugavet property, that is, every weakly compact polynomial $P\colon E \longrightarrow E$ satisfies the Daugavet equation $\|\id_{E} + P\| = 1+\|P\|$. The analogous conclusion also holds for the alternative Daugavet property.
	\end{abstract}
	
	\maketitle
	\thispagestyle{empty}
	
	\section{Introduction}
	A Banach space $X$ satisfies the \emph{Daugavet property} (\emph{DPr} in short) if
	the norm equality
	\begin{equation}\label{eq DP} \tag{DE}
		\|\id_{X} + T\| = 1+\|T\|
	\end{equation}
	holds for all rank-one bounded linear operators $T\colon X \longrightarrow X$ (equivalently, for all weakly compact linear operators on $X$, \cite{kssw}). Basic examples of Banach spaces with the Daugavet property include $C(K)$ spaces when $K$ is perfect, $L_1(\mu)$ spaces when $\mu$ is atomless, uniform algebras whose Choquet boundary is perfect, and isometric preduals of $L_1(\mu)$ spaces for which the set of extreme point of the dual ball is weak-star perfect up to rotations, among many others. We refer the interested reader to the seminal paper \cite{kssw} and to the more recent \cite{HallerLangemetsLimaNadelRueda, LopezRueda2021,MarRueda2022,RuedaTradaceteVillanueva} and references therein for more information and background.
	
	A very related property to the DPr is the following one. A Banach space $X$ satisfies the \emph{alternative Daugavet property} (\emph{ADP} in short) if the norm equality
	\begin{equation}\label{eq ADP} \tag{aDE}
		\max_{|w|=1} \|\id_{X} + w T\| = 1+\|T\|
	\end{equation}
	holds for all rank-one bounded linear operators $T\colon X \longrightarrow X$ (equivalently, for all weakly compact linear operators on $X$, \cite[Theorem~2.2]{MarOikh2004}). Here, the basic examples are $C(K)$ and $L_1(\mu)$ for all compact spaces $K$ and positive measures $\mu$, all uniform algebras, and all isometric preduals of $L_1(\mu)$, among many others. This property was formally introduced and deeply studied in \cite{MarOikh2004} (some related ideas had appeared before). We also refer the reader to \cite{AviKadMarMerSheTAMS2010,SpearsBook,Mar2008} for instance, for more information and background and for the relationship between the ADP and the study of the numerical index of Banach spaces.
	
	As this manuscript will deal with C$^*$-algebras and JB$^*$-triples, it makes sense to present the characterizations of the DPr and the ADP for these spaces given in \cite{BeMar2005,Mar2008,MarOikh2004,Oik2002}. We refer to Section~\ref{subsect:JBstar} for the definition of the involved concepts.
	
	\begin{theorem}[\textrm{\cite{BeMar2005,Mar2008,MarOikh2004,Oik2002}}]$ $ \newline \label{theorem:DPr-ADP-CstarJBstar}
		\textup{(a).} Let $X$ be a C$^*$-algebra.
		\begin{itemize}
			\item[(a1)] $X$ has the DPr if and only if $X$ is diffuse (i.e., it contains no atomic projections).
			\item[(a2)] $X$ has the ADP if and only if every atomic projection is central.
		\end{itemize}
		\textup{(b).} Let $X$ be a JB$^*$-triple.
		\begin{itemize}
			\item[(b1)] $X$ has the DPr if and only if $X$ contains no minimal tripotent.
			\item[(b2)] $X$ has the ADP if and only if every minimal tripotent is diagonalizing.
		\end{itemize}
	\end{theorem}
	
	This manuscript is devoted to study the natural extensions of the DPr and the ADP for polynomials. Let us introduce some notation. Even though both the DPr and the ADP have sense for both real and complex spaces, we will only deal with complex Banach spaces in this paper, as our main interest will be (complex) C$^*$-algebras and JB$^*$-triples. The closed unit ball of a Banach space $X$ will be denoted by $\mathcal{B}_{X}$. Let $X$ and $Y$ be Banach spaces. A (continuous) \emph{$N$-homogeneous polynomial} $P$ from $X$ to $Y$ is a mapping $P\colon X \longrightarrow Y$ for which we can find an operator $T\colon \overbrace{X\times \ldots \times X}^{N} \longrightarrow Y$ (continuous) multilinear and symmetric (i.e., $T(x_1,\ldots, x_N) = T(x_{\sigma(1)},\ldots, x_{\sigma(N)})$ for every permutation $\sigma$ of the set $\{1,\ldots, N\}$) satisfying $P(x) = T(x , \ldots , x)$ for every  $x \in X$. According to the usual notation, the symbol $\mathcal{P}(^N X,Y)$ will stand for the space of continuous $N$-homogeneous polynomials from $X$ to $Y$. Given $P\in \mathcal{P}(^N X,Y)$, we shall write $\widehat{P}$ for the unique continuous symmetric $N$-linear mapping associated with $P$. The scalar $N$-homogeneous polynomials on $X$ (i.e., $\mathcal{P}(^N X,\mathbb{C})$) will be simply denoted by $\mathcal{P}(^N X)$. A (general) polynomial from $X$ to $Y$ is a mapping $P\colon X\longrightarrow Y$ which can be written as a finite sum of homogeneous polynomials. We shall write $\mathcal{P}(X,Y)$ for the space of all polynomials from $X$ to $Y$, and $\mathcal{P}(X)$ for $\mathcal{P}(X,\CC)$. We say that $P\in \mathcal{P}(X,Y)$ is weakly compact when the closure of $P(\mathcal{B}_X)$ is weakly compact. We will consider $\mathcal{P}(X,Y)$ endowed with the usual supremum norm
	$$
	\| P \| = \sup\{ \| P(x)\| \colon \|x\| \leq 1\}.
	$$
	Using this norm, it makes sense to consider the equations \eqref{eq DP} and \eqref{eq ADP} in the space $\mathcal{P}(X,X)$ and hence to consider, as it is done in \cite{ChoGarMaestreMartin2007,cgmm2,Santos-studia}, the following two properties. We say that $X$ has the \emph{polynomial Daugavet property} \cite{ChoGarMaestreMartin2007,cgmm2} if every weakly compact $P\in \mathcal{P}(X,X)$ satisfies \eqref{eq DP}. $X$ has the \emph{alternative polynomial Daugavet property} if every weakly compact $P\in \mathcal{P}(X,X)$ satisfies \eqref{eq ADP}. It is immediate that the polynomial versions of the Daugavet and the alternative Daugavet properties imply the usual ones, respectively.
	
	It is an open problem whether the Daugavet property implies the polynomial Daugavet property, but it is known that this is the case in some families of Banach spaces as $C(K)$ spaces and some generalizations \cite{ChoGarMaestreMartin2007,cgmm2}, for $L_1(\mu)$ spaces and for vector-valued $L_1$ spaces \cite{mmp}, and also for isometric preduals of $L_1(\mu)$, for uniform algebras and for spaces of Lipschitz functions \cite{MarRueda2022}, among others. For the case of the ADP, the situation is different as, for instance, the complex space $\ell_1$ fails the alternative polynomial Daugavet property \cite{ChoGarMaestreMartin2007} (despite of the fact that $\ell_1$ has the ADP). It is however known that every complex $C(K)$ space has the alternative polynomial Daugavet property \cite{ChoGarMaestreMartin2007}. The real case is even worse: real $c_0$ and real $\ell_1$ fail the alternative polynomial Daugavet property \cite{ChoGarMaestreMartin2007}. E.~Santos established in \cite{Santos2014} that every continuous polynomial of finite type or approximable on a JB$^*$-triple having the DPr (respectively, the ADP) satisfies the Daugavet equation \eqref{eq DP} (respectively, the alternative Daugavet equation \eqref{eq ADP}).
		
	The main aim of this paper is to show that for C$^*$-algebras and JB$^*$-triples, the DPr implies the polynomial Daugavet property and the ADP implies the alternative polynomial Daugavet property.
	
\begin{theorem}[Main result] Let $X$ be a JB$^*$-triple {\rm(}in particular, a C$^*$-algebra{\rm)}.
\begin{enumerate}[$(a)$]
\item If $X$ has the DPr, then it has the polynomial Daugavet property.
\item If $X$ has the ADP, then it has the alternative polynomial Daugavet property.
\end{enumerate}
\end{theorem}
	
It should be noticed that this result provides a complete positive solution to the problems posed by E.~Santos in \cite[\S 4]{Santos2014}.
	
The proofs of these two statements will appear in Theorem~\ref{t DPr for general wk polynomials} and Theorem~\ref{t ADPr for general wk polynomials}, respectively, where we actually rediscover a new proof of  Theorem~\ref{theorem:DPr-ADP-CstarJBstar}.
	
Let us comment that the main idea to get the polynomial Daugavet property in $C(K)$ spaces and in some related classes of Banach spaces is to produce a number of ``good'' $c_0$-sequences in the space (or in the bidual), and then use that continuous scalar polynomials in $c_0$ are weakly continuous on bounded sets. This idea does not seem to work for noncommutative C$^*$-algebras nor for general JB$^*$-triples. Therefore, we need to find a suitable substitute of it, and the alternative tool is the sequential continuity of scalar polynomials on JB$^*$-triples (and so on C$^*$-algebras) for the strong$^*$-topology (see Corollary~\ref{c scalar polynomials on JBstar triples as sequentilly strong* continuous}), a result which is of interest by itself. The remaining strategy to show that JB$^*$-triples having the DPr and the ADP actually satisfy their polynomial versions consists in extending and improving ideas from the original proofs given in \cite{MarOikh2004,Oik2002} for the linear properties. There is however a substantial turn here, we do not rely the arguments in controlling the norm of a determined finite set of linear functionals when restricted to the Peirce-2 and Peirce-1 subspaces of an appropriate tripotent in the bidual, instead of that we shall show that  the joint strong$^*$-continuity of the triple product in the bidual spaces simplifies the computation of the norms for scalar polynomials via pointwise convergence.
	
The outline of the paper is as follows. Section \ref{subsect:polDPR-polADP} is devoted to revisit some known tools and results to study the Daugavet equation for polynomials on a Banach space $X$ which are weakly continuous on bounded sets. We prove that, assuming that $X$ satisfies the DPr, every weakly compact polynomial on $X$ which is weakly continuous on bounded sets satisfies the Daugavet equation (see Theorem \ref{prop:generalDPr-weaklycontinouspolynomials}). Section \ref{subsect:JBstar} contains a very brief introduction to JB$^*$-triples together with some basic tools and results required in our arguments (like the strong$^*$-topology and its main properties). We shall also revisit the results guaranteeing that every scalar polynomial on a general JB$^*$-triple is sequentially strong$^*$-continuous. Finally, section~\ref{section:proofs} is devoted to present the proofs of the main results.

	\section{Revisiting the Daugavet equation for polynomials}\label{subsect:polDPR-polADP}
	Our goal here is to recall some basic facts on the polynomial Daugavet property and the alternative polynomial Daugavet property which will be useful in the sequel. We start with characterizations of both properties in terms of polynomials of rank-one.
	
	We start with a result for the polynomial Daugavet property.
	\begin{proposition}[\textrm{\cite[Theorem 1.1 or Corollary 2.2]{cgmm2}}]\label{prop_charpolDP}
		Let $X$ be a Banach space. Then, the following are equivalent:
		\begin{enumerate}[$(a)$]
			\item $X$ has the polynomial Daugavet property {\rm(}i.e., every weakly compact polynomial on $X$ satisfies \eqref{eq DP}{\rm)}.
			\item Every $P\in \mathcal{P}(X,X)$ of the form $x\longmapsto p(x)a$ for suitable $p\in \mathcal{P}(X)$ and $a\in X$, satisfies \eqref{eq DP}.
			\item Given $p\in \mathcal{P}(X)$ with $\|p\|=1$, $a\in X$ with $\|a\|=1$, and $\varepsilon>0$, there are $x\in \mathcal{B}_X$ and $w\in \CC$ with $|w|=1$ such that
			$$
			\re w p(y)>1-\varepsilon \quad \text{and} \quad \|a+ w x\|>2-\varepsilon.
			$$
		\end{enumerate}
	\end{proposition}
	
	The analogous result for the alternative polynomial Daugavet property also holds.
	\begin{proposition}[\textrm{\cite[Corollary~1.2]{cgmm2}}]\label{prop_charpol-ADP}
		Let $X$ be a Banach space. Then, the following are equivalent:
		\begin{enumerate}[$(a)$]
			\item $X$ has the alternative polynomial Daugavet property {\rm(}i.e., every weakly compact polynomial on $X$ satisfies \eqref{eq ADP}{\rm)}.
			\item Every $P\in \mathcal{P}(X,X)$ of the form $x\longmapsto p(x)a$ for suitable $p\in \mathcal{P}(X)$ and $a\in X$, satisfies \eqref{eq ADP}.
			\item Given $p\in \mathcal{P}(X)$ with $\|p\|=1$, $a\in X$ with $\|a\|=1$, and $\varepsilon>0$, there are $x\in \mathcal{B}_X$ and $w\in \CC$ with $|w|=1$ such that
			$$
			|p(x)|>1-\varepsilon \quad \text{and} \quad \|a+ w x\|>2-\varepsilon.
			$$
		\end{enumerate}
	\end{proposition}
	
	As we already mentioned in the introduction, an strategy to prove the polynomial Daugavet property in some Banach spaces like $C(K)$ spaces,  uniform algebras, or isometric preduals of $L_1$-spaces is to show that these spaces (or their biduals) are full of ``good copies'' of $c_0$, and then use the fact that scalar polynomials on $c_0$ are weakly continuous on bounded sets (see \cite[Proposition~6.3]{cgmm2} and \cite[Proposition~4.3]{MarRueda2022}). Of course this strategy cannot be used in all Banach spaces with the DPr (for instance, $L_1[0,1]$ contains no copies of $c_0$). Here we can get a partial result showing that in a Banach space with the DPr, weakly compact polynomials which are weakly continuous on bounded sets satisfy \eqref{eq DP}. As far as we know, this result is new.
	
	\begin{theorem}\label{prop:generalDPr-weaklycontinouspolynomials}
		Let $X$ be a Banach space with the DPr. Then, given $p\in \mathcal{P}(X)$ which is weakly continuous on bounded sets and satisfies $\|p\|=1$, $a\in X$ with $\|a\|=1$, and $\varepsilon>0$, there are $x\in \mathcal{B}_X$ and $w\in \CC$ with $|w|=1$ such that
		$$
		\re w p(x)>1-\varepsilon \quad \text{and} \quad \|a+ w x\|>2-\varepsilon.
		$$
		As a consequence, every weakly compact $P\in \mathcal{P}(X,X)$ which is weakly continuous on bounded sets satisfies \eqref{eq DP}.
	\end{theorem}
	
	We will make use of \cite[Lemma~3]{ShvydkoyJFA2000} for which only the following weaker version will be needed.
	
	\begin{lemma}[\textrm{\cite[Lemma~3]{ShvydkoyJFA2000}}]\label{lemma-Shvidkoy}
		Let $X$ be a Banach space with the DPr. Then, for every $y\in \mathcal{B}_X$ with $\|y\|=1$ and every $\varepsilon>0$, the set
		$$
		\{x\in \mathcal{B}_X\colon \|y+x\|>2-\varepsilon\}
		$$
		is weakly dense in $\mathcal{B}_X$.
	\end{lemma}
	
	Let us observe that the lemma actually provides a characterization of Banach spaces with the Daugavet property, since the reciprocal implication is immediate by the well known characterization of the DPr using slices (see \cite[Lemma~2.1]{kssw}).
	
	\begin{proof}[Proof of Theorem~\ref{prop:generalDPr-weaklycontinouspolynomials}]
		Take $w\in \CC$ with $|w|=1$ such that $\sup\limits_{x\in \mathcal{B}_X} \re w p(x)>1-\varepsilon$ and observe that, being $p$ weakly continuous on $\mathcal{B}_X$, the set
		$$
		U:=\{x\in \mathcal{B}_X\colon \re w p(x)>1-\varepsilon\}
		$$
		is non-empty and weakly open relative to $\mathcal{B}_X$. Now, Lemma~\ref{lemma-Shvidkoy}, applied to $\bar{w}a,$ assures that there is $x_0\in U$ such that $\bigl\|\bar{w}a + x_0\bigr\|>2-\varepsilon$. If clearly follows that
		\begin{equation*}
			\re w p(x_0)>1-\varepsilon  \quad \text{and} \quad \|a+ w x_0\|>2-\varepsilon.
		\end{equation*}
		To prove the last statement in the theorem, we just have to follow the lines of the proof of $(c)\Rightarrow (a)$ in Proposition~\ref{prop_charpolDP} given in \cite[Theorem~1.1]{cgmm2}, and observe that in order to prove that every weakly compact and weakly continuous on bounded sets  $P\in \mathcal{P}(X,X)$ satisfies \eqref{eq DP}, we only need to check that $(c)$ holds for scalar polynomials which are weakly continuous on bounded sets, that is, we only need what we already proved in the first part of this proposition.
	\end{proof}
	
	Some remarks are worth mentioning.
	
	\begin{remarks}$ $
		\begin{enumerate}[$\bullet$]
			\item Theorem~\ref{prop:generalDPr-weaklycontinouspolynomials} cannot be applied to show that C$^*$-algebras or JB$^*$-triples with the DPr have the polynomial Daugavet property. Actually, it cannot be used to get the polynomial Daugavet property from the DPr in any Banach space. Indeed, every Banach space with the DPr contains an isomorphic copy of $\ell_1$ \cite[Theorem~2.9]{kssw} and it is known that every Banach space containing $\ell_1$ admits a continuous scalar polynomial which is not weakly continuous on bounded sets {\rm (}cf.\ \cite[Proposition~2.36 and its proof]{Di99}{\rm )}.
			\item In case that every continuous scalar polynomial on $X$ is weakly sequentially continuous  {\rm (}as happens in Banach spaces with the Dunford-Pettis property, see \cite[Proposition~2.34]{Di99}{\rm )}, then a ``sequential'' version of Lemma~\ref{lemma-Shvidkoy} for such an $X$ {\rm (}i.e., weak sequential denseness instead of weak denseness{\rm )} would be enough to get the polynomial Daugavet property from the DPr on $X$. We do not know if such a sequential version of Shvidkoy lemma is true for C$^*$-algebras or JB$^*$-triples, but it is certainly not true for all Banach spaces with the DPr: there is a Banach space satisfying both the DPr and the Schur property \cite{KW-shurDPr}.
			\item It was proved in \cite[Theorem~3.5]{Santos2014} that when $X$ is a C$^*$-algebra or a JB$^*$-triple and $X$ has the DPr, then the polynomials of finite-type on $X$ satisfy \eqref{eq DP}. Recall that a polynomial $P\in \mathcal{P}(^N X,Y)$ is said to be of \emph{finite type} if there exists a finite subset $\{\varphi_i\}_{i=1}^{m}\subset X^*$ and elements $y_i$'s in $Y$ such that $\displaystyle P(x) = \sum_{i=1}^m y_i \varphi_i (x)^N$ {\rm (}$x\in X${\rm )}. It is clear that polynomials of finite type are weakly continuous {\rm(}in particular, weakly continuous on bounded sets{\rm)}, so the cited result \cite[Theorem~3.5]{Santos2014} actually follows from our Theorem~\ref{prop:generalDPr-weaklycontinouspolynomials} and it is indeed true for all Banach spaces with the DPr, not only for C$^*$-algebras and JB$^*$-triples.
			\item The analogous result to Theorem~\ref{prop:generalDPr-weaklycontinouspolynomials} for the ADP and the polynomial alternative Daugavet property is not true. Indeed, the {\rm (}real or complex{\rm )} two-dimensional space $\ell_1^2$ has the ADP \cite{MarOikh2004}, but fails the polynomial alternative Daugavet property {\rm (}see Example 3.12 and the paragraph below Remark 3.13 in \cite{ChoGarMaestreMartin2007}{\rm )}. Being finite-dimensional, it is immediate that all polynomials on $\ell_1^2$ are weakly continuous {\rm(}on bounded sets{\rm)}. By the way, this shows that the result analogous to Shvidkoy's Lemma \ref{lemma-Shvidkoy} for the ADP cannot be true.
		\end{enumerate}
	\end{remarks}
	
	\section{A closer look at the theory of C$^*$-algebras and JB$^*$-triples}\label{subsect:JBstar}
	
	This section is not only devoted to summarize the basic definitions of JB$^*$-triples, but to revisit and adapt the geometric and topological tools required in our arguments. We shall also rediscover and extend some geometric tools employed in previous references from a quite different point of view.
	
	There are several worthy arguments to consider and study the complex Banach spaces in the category of JB$^*$-triples. First, the class of JB$^*$-triples is strictly wider than the class of C$^*$-algebras, and in many cases the geometric properties are easier affordable from this more general point of view. Perhaps the most important motivation to introduce JB$^*$-triples arises in holomorhic theory and the classification of bounded symmetric domains in arbitrary complex Banach spaces, domains that play the role of simply connected domains in Riemann mapping theorem (cf.\ \cite{Ka83}). Let us briefly recall the definition. A complex Banach space $E$ is called a \emph{JB$^*$-triple} if it admits a continuous triple product $\J \cdot\cdot\cdot\colon
	E\times E\times E \longrightarrow E,$ which is symmetric and bilinear in the first and third variables, conjugate-linear in the middle one,
	and satisfies the following axioms:
	\begin{enumerate}[{\rm (a)}] \item (Jordan identity)
		$$L(a,b) L(x,y) = L(x,y) L(a,b) + L(L(a,b)x,y)
		- L(x,L(b,a)y)$$ for $a,b,x,y$ in $E$, where $L(a,b)$ is the operator on $E$ given by $x \longmapsto \J abx;$
		\item $L(a,a)$ is a hermitian operator with non-negative spectrum for all $a\in E$;
		\item $\|\{a,a,a\}\| = \|a\|^3$ for each $a\in E$.
\end{enumerate}
	
The triple product \begin{equation}\label{eq C* triple product} \{a,b,c\} = \frac12 ( ab^* c + c b^* a)
\end{equation} provides a structure of JB$^*$-triple for every C$^*$-algebra and every closed subspace of $B(H)$ which is closed for the triple product just commented -- in particular the space $B(H,K)$ of all bounded linear operators between two complex Hilbert spaces and all complex Hilbert spaces are JB$^*$-triples. There exist JB$^*$-triple which cannot be embedded as JB$^*$-subtriples of $B(H)$, they are related to the so-called exceptional Cartan factors (cf.\ \cite{FriRu86}). A JBW$^*$-triple is a JB$^*$-triple which is also a dual Banach space. In analogy with Sakai's theorem, every JBW$^*$-triple admits a unique (isometric) predual and its triple product is separately weak$^*$ continuous (cf.\ \cite{BarTi}).
	
Let $A$ be a C$^*$-algebra. It is known that the fixed points of the triple product \eqref{eq C* triple product} are precisely the partial isometries in $A$ (i.e., those $e\in A$ such $e e^* e = e$, equivalently, $ee^*$ or $e^* e$ is a projection). If we fix a partial isometry $e\in A,$ we can easily decompose $A$ in the form $$A = ee^* A e^*e + \left( (1-ee^*) A e^*e + ee^* A (1-e^*e) \right) + (1-ee^*) A (1-e^*e), $$ decomposition which is known under the name of Peirce decomposition of $A$ associated to $e$. Clearly, projections in $A$ are nothing but those positive partial isometries in $A$.

More generally, if we have a JB$^*$-triple $E$, the elements $e\in E$ satisfying $\{e,e,e\}=e$ are called \emph{tripotents}.  Each tripotent $e$ in $E$ produces a \emph{Peirce decomposition} of the space $E$ in terms of the eigenspaces of the operator $L(e,e)$:
	\begin{equation}\label{Peirce decomp} {E} = {E}_{2} (e) \oplus  {E}_{1} (e) \oplus {E}_{0} (e),\end{equation} where ${
		E}_{k} (e) := \{ x\in {E}\colon L(e,e)x = {\frac k 2} x\}$ is a subtriple of ${E}$ called the \emph{Peirce-$k$ subspace} ($k=0,1,2$). The natural projection of ${E}$ onto ${E}_{k} (e)$ is known as the  \emph{Peirce-$k$ projection}, and it is usually denoted by $P_{k} (e)$. We shall later employ that Peirce projections are all contractive (cf.\ \cite[Corollary 1.2]{FriRu85}).
	
A tripotent $e$ in a JB$^*$-triple $E$ is called \emph{minimal} if $e$ is non-zero and $E_2(e) = \mathbb{C} e.$
	
Triple products among elements in different Peirce subspaces obey certain laws known as \emph{Peirce arithmetic}. Concretely, the inclusion $\J {{E}_{k}(e)}{{E}_{l}(e)}{{E}_{m}(e)}\! \subseteq {E}_{k-l+m} (e),$ and the identity $\J {{E}_{0}(e)}{{E}_{2} (e)}{{E}}\! =\! \J {{E}_{2} (e)}{{E}_{0} (e)}{{E}}\! =\! \{0\},$ hold for all $k,l,m\in \{0,1,2\}$, where ${E}_{k-l+m} (e) = \{0\}$ whenever $k-l+m$ is not in $\{0,1,2\}$. The Peirce-$2$ subspace ${E}_{2} (e)$ is a unital JB$^*$-algebra with respect to the product and involution given by $x \circ_e y = \J xey$ and $x^{*_e} = \J exe,$ respectively (cf.\ \cite[Corollary~4.2.30]{Cabrera-Rodriguez-vol1}).
	
	A tripotent $e$ in a JB$^*$-triple $E$ is said to be \emph{diagonalizing} if $E_1 (e) =\{0\}$. Let us observe that a projection $p$ in a C$^*$-algebra $A$ is diagonalizing as tripotent if and only if it is a central projection.
	
	Projections $p,q$ in a C$^*$-algebra $A$ are called orthogonal if $p q =0.$ If we consider general elements $a,b\in A,$ according to the usual notation, we shall say that $a$ and $b$ are \emph{orthogonal} ($a\perp b$ in short)  if $a b^* = b^* a= 0$. It is known that $a\perp b$ if and only if $L(a,b)=0$ (cf.\ \cite[Lemma 1]{BurFerGarMarPe08}). In a general JB$^*$-triple $E$, elements $a,b$ are called \emph{orthogonal} if $L(a,b) =0$. In case that $a$ and $b$ are tripotents, it is easy to see that $a\perp b$ if and only if $a\in E_0(b)$, if and only if $b\in E_0(a)$. The second identity in the so-called Peirce arithmetic precisely tells that $a\perp b$ whenever $a\in E_0(e)$ and $b\in E_2(e)$ for a tripotent $e\in E.$ Another geometric property of orthogonal elements in JB$^*$-triples assures that every two of orthogonal elements $a,b$ in a JB$^*$-triple are $M$-orthogonal, that is, \begin{equation}\label{eq M-orhtogonality} \hbox{$\|a\pm b\| = \max\{\|a\|,\|b\|\}$ (cf.\ \cite[Lemma 1.3$(a)$]{FriRu85}).}
	\end{equation}
	
	There is a natural partial ordering among tripotents elements in a JB$^*$-triple $E$ given by $e\leq u$ if $u-e$ is a tripotent and $(u-e)\perp e.$
	
	The Gelfand-Naimark-Segal construction is one of the best known results in representation theory of C$^*$-algebras. One of the many ideas in this construction provides a tool to define Hilbertian structures associated with positive functionals in the dual of a C$^*$-algebra. Concretely, if $\phi$ is a positive functional in the dual of a C$^*$-algebra $A$, the mapping $(a,b)\longmapsto \phi (b^* a)$ defines a semi-positive sesquilinear form on $A$. The symmetric version of these sesquilinear form (i.e., $(a,b)\longmapsto \phi (\frac{b^* a + a^* b}{2})$) defines the pre-Hilbertian seminorm of the form $\|a\|_{\phi}^2 = \phi (\frac{a^* a + a a^*}{2})$ ($a\in A$) appearing in the non-commutative version of Grothendieck's inequalities by G.~Pisier and U.~Haagerup (cf.\ \cite{Pi78, Haa85, Pi2012}).
	
	Despite of the lacking of a positive cone in general JB$^*$-triples, J.~T.~Barton and Y.~Friedman established in \cite[Proposition 1.2]{BarFri1987Grothendieck} the following procedure to define a pre-Hilbertian seminorm associated with a functional $\varphi$ in the dual, $E^*$, of a JB$^*$-triple $E$: for each $z\in E^{**}$ with $\|z\| =1$ and $\varphi (z)=\|\varphi\|,$ the mapping $$\begin{aligned} E\times & E \xrightarrow{\hspace*{0.59161cm}}  \mathbb{C} \\
		(x,y)&\longmapsto {\varphi} \J xyz
	\end{aligned}$$ is a semi-positive sesquilinear form on $E,$ which does not depend on the choice of the element $z\in E^{**}$. The corresponding prehilbertian seminorm on $E$ is denoted by  $\|x\|_{{\varphi}}^2 := \varphi \{x,x,z\}$ $(x\in E)$. All pre-Hilbertian seminorms associated with positive functionals in the dual space of a C$^*$-algebra $A$ arise in this way, since for each positive functional $\phi$ we have $\phi(\mathbf{1})=\|\phi\|$ for the unit in $A^{**}$, and hence $\|x\|_{\phi}^2 = \phi \{x,x,\mathbf{1}\} = \phi (\frac{x^* x + x x^* }{2})$. These pre-Hilbertian seminorms play a fundamental role in the results known as Grothedieck's inequalities for JB$^*$-triples and in the definition of the strong$^*$-topology which will be recalled later (cf.\ \cite{BarFri90, HamKalPePfi20Groth}).
	
	It is known that
	$$|\varphi (x)| \leq \|\varphi\| \ \|x\|_{\varphi}$$ for all $\varphi\in E^*$, $x\in E$ (see \cite[comments before Definition 3.1]{BarFri90}).
	
	Among the results we shall revisit here, we shall present a quantitative version of results stated by M.~Mart{\'i}n and T.~Oikhberg \cite[Lemma 4.16]{MarOikh2004} and E.~Santos \cite[Lemma 3.3]{Santos2014}. We begin by recalling a result borrowed from \cite{FerPe06}, which offers a technical control of a functional on Peirce-1 and -2 subspaces associated with a tripotent at which the pre-Hilbertian seminorms of the functionals are ``small''.
	
	\begin{lemma}\label{l tecnical seminorm}{\rm\cite[Lemma 3.2]{FerPe06}} Let $e$ be a tripotent in a JB$^*$-triple $E$, and let $\varphi$ be a norm-one element in $E^{*}$ such that $\|e\|_{\varphi}^2 < \delta$. Then the inequalities $$ |\varphi P_2 (e) (x)| < 3 \sqrt{\delta} \ \|P_2 (e)
		(x)\|, \hbox{ and } |\varphi P_1 (e) (x)| < 6 \sqrt{\delta} \ \|P_1 (e)
		(x)\|,$$ hold for all $x\in E$.
	\end{lemma}
	
	It is known, and easily deducible from the definition of orthogonality, that the square of each pre-Hilbertian seminorm behaves additively on orthogonal elements, that is, if $a_1,\ldots, a_n$ are mutually orthogonal elements in a JB$^*$-triple $E$ and $\varphi\in E^*$, we have $\displaystyle\Bigl\| \sum_{j=1}^n a_j \Bigr\|_{\varphi}^2 =  \sum_{j=1}^n \left\|a_j\right\|_{\varphi}^2.$ Our next goal is a simplified proof of \cite[Lemma~4.16]{MarOikh2004} and \cite[Lemma~3.3]{Santos2014}, which additionally provides a quantitative conclusion. The desired statement is analyzed in the following remark.

	\begin{remark}\label{r simplified statement and proof of OM and Sant} For each $\varepsilon >0$ and $k\in \mathbb{N}$ take a natural number $n_0$ satisfying $n_0 > \frac{12^2 k}{\varepsilon^2}$. Then for every finite family $\{\varphi_1,\ldots, \varphi_k\}$ of non-zero functionals in the closed unit ball of the dual space of a JB$^*$-triple $E,$ and each finite family $\{e_1,\ldots, e_M\}$ of non-zero mutually orthogonal tripotents in $E$ with $M\geq n_0$, there is $m\in \{1,\ldots, M\}$ satisfying $\|\varphi_i|_{E_2(e_m)\oplus E_1(e_m)}\| = \|\varphi_i (P_2(e_m)+ P_1(e_m))\|< \varepsilon,$ for all $1\leq i \leq k$.
		
		Namely, let us fix $\varepsilon>0,$ $k\in \mathbb{N}$. Find $n_0$ in $\mathbb{N}$ such that $\frac{k}{n_0}<\frac{\varepsilon^2}{12^2}.$ Consider functionals and mutually orthogonal tripotents as in the statement with $M\geq n_0$, and define the following pre-Hilbertian seminorm:
		$$
		\|x\|_{\varphi_1,\ldots, \varphi_k}^2 := \sum_{j=1}^k \|x\|_{\varphi_j}^2, \ \ (x\in E).
		$$
		Since the square of every $\|\cdot\|_{\varphi_j}$ is additive on orthogonal elements, the square of the seminorm $\|\cdot\|_{\varphi_1,\ldots, \varphi_k}$ is also additive on orthogonal elements. Having in mind that the element $\sum_{i=1}^M e_i$ is a tripotent and the fact that $e_1,\ldots, e_M$ are mutually orthogonal, we deduce that
		$$ \sum_{i=1}^M  \left\| e_i \right\|_{\varphi_1,\ldots, \varphi_k}^2 = \Bigl\| \sum_{i=1}^M e_i \Bigr\|^2_{\varphi_1,\ldots, \varphi_k} = \sum_{j=1}^k \Bigl\| \sum_{i=1}^M e_i \Bigr\|^{2}_{\varphi_j} \leq \sum_{j=1}^k \|\varphi_j\| \Bigl\| \sum_{i=1}^M e_i  \Bigr\|^2 \leq k .
		$$
		Therefore, there exists $m\in \{1,\ldots, M\}$ such that $\left\| e_m \right\|_{\varphi_1,\ldots, \varphi_k}^2 <\frac{k}{M} \leq \frac{k}{n_0} < \frac{\varepsilon^2}{12^2}.$ Since for each $i \in \{1, \ldots, k\}$ we have $\|e_m\|_{\varphi_i}^2 \leq \left\| e_m \right\|_{\varphi_1,\ldots, \varphi_k}^2< \frac{\varepsilon^2}{12^2},$ Lemma \ref{l tecnical seminorm} assures that \[\begin{aligned}|\varphi_i \left(P_2 (e_m) + P_1 (e_m)\right) (x)| &< 3 \frac{\varepsilon}{12} \ \left\|P_2 (e_m) (x)\right\|+ 6 \frac{\varepsilon}{12} \ \left\|P_1 (e_m)
			(x)\right\| \\
			&\leq 9 \frac{\varepsilon}{12} \ \left\|P_2 (e_m) (x) + P_1 (e_m) (x)\right\|
		\end{aligned}\] for all $1\leq i\leq k$ and all $x\in E$, where in the last inequality we applied that Peirce projections are contractive.
	\end{remark}

The arguments in previous studies, like \cite{MarOikh2004,Santos2014} rely on controlling the norms of the restriction to the Perice-2 and Peirce-1 subspaces associated with a tripotent in the bidual space of all elements a finite family of functionals in the dual space of a JB$^*$-triple $E,$ as in the previous remark. The same conclusion does not seem easy achievable for an arbitrary scalar polynomial. For the proof in this paper we shall rely on a weaker pointwise convergence.

One of the novelties in this note, compared with previous forerunners,  consists in replacing finite families of mutually orthogonal elements by sequences of mutually orthogonal elements. In order to sum countable families of mutually orthogonal tripotents, we need to employ an appropriate topology to assure the convergence. Actually, the following stronger property holds: every arbitrary family $\{e_i\}_{i\in \Lambda}$ of mutually orthogonal tripotents in a JBW$^*$-triple $M$ is summable with respect to the weak$^*$ topology of $M$, and its limit, denoted by $\displaystyle \sum_{i\in \Lambda} e_i = \hbox{w$^*$-}\sum_{i\in \Lambda} e_i$, is another tripotent in $M$ (cf.\ \cite[Corollary 3.13]{Horn87} or \cite[Proposition 3.8]{Batt91}).
	
	We shall also make use of the strong$^*$-topology of JB$^*$-triples. For this reason, it is worth to recall some basic facts on this topology. According to \cite[\S 3]{BarFri90}, the \emph{strong$^*$-topology} of a JB$^*$-triple $E$ is the topology determined by all the seminorms of the form $\|\cdot\|_{\varphi},$ with $\varphi$ running in the dual space of $E$ (or equivalently, in the unit sphere of $E^*$). If $M$ is a JBW$^*$-triple, the strong$^*$-topology of $M$ is the topology determined by all the seminorms of the form $\|\cdot\|_{\varphi}$ with $\varphi$ running in the predual space of $M$. Let us observe that the strong$^*$-topology of $E$ is precisely the restriction to $E$ of the strong$^*$-topology of $E^{**}$. This topology enjoys some useful properties, it is compatible with the duality $(E,E^*)$ and stronger than the weak topology of $E$.
	
	When a C$^*$-algebra $A$ is regarded as a JB$^*$-triple, the strong$^*$-topology on $A$ in the triple sense coincides with the strong$^*$ topology in the usual C$^*$- sense (cf.\ \cite[1.8.7]{Sak} and  \cite[p.\ 258-259]{BarFri90}). In what concerns this paper, we remark that the triple product is jointly strong$^*$-continuous on bounded sets (cf.\ \cite{Rod89}, \cite[\S 4]{PeRo01}). It is well known that every sequence of mutually orthogonal projections in a C$^*$-algebra is strong$^*$-null, and the same property holds for every sequence of mutually orthogonal tripotents in a JB$^*$-triple (cf.\ \cite[Comments in page 86]{PePfi15}).\label{eq orthogonal sequences are strong* null}
		
	A (closed) subtriple $I$ of a JB$^*$-triple $E$ is said to be an \emph{ideal} (respectively, an \emph{inner ideal}) of $E$ if $\J EEI + \J EIE\subseteq I$ (respectively, $\J IEI \subseteq I$). For example, if $p,q$ are two projections in a C$^*$-algebra $A$, the subtriple $p A q$ is an inner ideal of $A$, in particular, the Peirce-2 subspace associated with a partial isometry is an inner ideal.
	
	It is shown in \cite[Lemma 3.2]{Santos2014} that given two tripotents $e, e_1$ in a JB$^*$-triple $E$ the conditions $e\geq e_1$ (i.e., $e = e-e_1 = e_1$ with $e-e_1\perp e-1$) and $e_1$ minimal in the Peirce-2 subspace $E_2(e)$, imply that $e_1$ is minimal in $E$. This is actually a consequence of the following more general property, where we are not assuming any order relationship among the tripotents.
	
	\begin{lemma}\label{l minimal tripotents in inner ideals} Let $I$ be an inner ideal of a JB$^*$-triple $E$. Suppose that $e_1$ is a minimal tripotent in $I$. Then $e_1$ is minimal in $E$. In particular, if $e_1$ is a minimal tripotent in $E_2(e)$ or in $E_0(e)$ for a fixed tripotent $e\in E$, then $e_1$ is a minimal tripotent in $E$.
	\end{lemma}
	
	\begin{proof} By assumptions $I_2(e_1)= \mathbb{C} e_1$. Each element $x$ in the unital JB$^*$-algebra $E_2(e_1)$ decomposes in the form $x = h + i k$ with $h,k$ self-adjoint elements in $E_2(e_1)$. In particular, $k = k^{*_{e_1}} = \{e_1,k, e_1\}$ and $h = h^{*_{e_1}} = \{e_1, h, e_1\}$. Having in mind that $I$ is an inner ideal with $e_1\in I$, we get $k = \{e_1,k, e_1\},$ $h = \{e_1, h, e_1\}\in I$. Thus $h$ and $k$ are self-adjoint elements in the JB$^*$-algebra $I_2 (e_1) = \mathbb{C} e_1,$ and consequently $x = h +i k \in \mathbb{C} e_1$, which gives the desired statement. The rest is clear because, for each tripotent $e\in E$, the Peirce subspaces $E_2(e)$ and $E_0(e)$ are inner ideals.
	\end{proof}
	
The C$^*$-subalgebra of C$^*$-algebra generated by a single non-normal element is not always representable as a C$^*$-algebra of continuous functions vanishing at infinity on a locally compact Hausdorff space. However, the JB$^*$-subtriple, $E_a$, generated by a single element $a$ in a JB$^*$-triple $E$ is JB$^*$-triple isomorphic to $C_0(Sp(a))$ for some unique compact set $Sp(a) \subseteq [0,\|a\|]$ with $\|a\| \in Sp(a)$, where $C_0(Sp(a))$ denotes the commutative C$^*$-algebra of all continuous complex-valued functions on $Sp(a)$ vanishing at zero if $0\in Sp(a)$. It is further known that we can actually find a triple isomorphism $\Psi_a\colon E_a \longrightarrow C_0(Sp(a))$ mapping $a$ to the natural inclusion of $Sp(a)$ into $\mathbb{C}$ (cf.\ \cite[Corollary 1.15]{Ka83} and \cite[Proposition 3.5$(iii)$]{Ka96}). The set $Sp(a)$ is called the  \emph{triple spectrum of $a$}. It is worth to note that for $a\neq 0$, the triple spectrum $Sp(a)$ is precisely the set of all positive square roots of the elements in the spectrum of the operator $L(a,a)|_{E_a}$ in the unital complex Banach algebra $B(E_a)$, that is, $$\begin{aligned}
	Sp(a) &=\{t\in \mathbb{C}\colon t^2 \in \sigma_{B(E_a)} (L(a,a)|_{E_a}) \} \\
	&= \{ t \in \mathbb{C} \colon L(a,a)|_{E_a}- t^2 Id_{E_a} \hbox{ not invertible in }  B(E_a) \}\\
	&= \{ t \in \mathbb{C}\colon a \notin (L(a,a)-t^2 Id) (E)\}
\end{aligned}$$ 
(cf.\ \cite[Corollary 3.4]{Ka96}). For $a =0$ we set $Sp(a) = \{0\}$ --in such a case $B(E_a)$ is not a unital Banach algebra.
	
Orthogonality is the key notion to define the rank of a JB$^*$-triple $E$. A subset $\mathcal{S}\subset E$ is \emph{orthogonal} if $0\notin \mathcal{S}$ and $a\perp b$ for all $a,b\in \mathcal{S}$. The minimal cardinal number $r$ satisfying $\hbox{card}(S) \leq r$ for every orthogonal subset $S \subseteq E$ is called the \emph{rank} of $E$ ($r(E)$ in short). The rank of a tripotent $e$ in $E$, $r(e)$, is defined as the rank of the Peirce-2 space $E_2(e)$. A JB$^*$-triple has finite rank if and only if it is reflexive (cf.\ \cite[Proposition 4.5]{BuChu92} and \cite[Theorem 6]{ChuIo90} or \cite{BeLoRo03, BeLoPeRo04}).
	
	We would like to remark some knwon facts on tripotents of finite rank. We observe that every minimal tripotent in a JBW$^*$-triple $M$ lies in the atomic part of $M$, actually the atomic part of $M$ is the weak$^*$-closed ideal of $M$ generated by all minimal tripotents in $M$ (cf.\ \cite{FriRu85}). The atomic part of $M$ coincides with an $\ell_{\infty}$-sum of Cartan factors (see \cite[Proposition 2]{FriRu85} and \cite[Corollary 1.8]{Horn87b}). By working on each of the summands appearing in the atomic part of $M$, we deduce that if a tripotent $e$ writes as the orthogonal sum of a finite family of $n$ minimal tripotents, then it has finite rank $n,$ and reciprocally, if a tripotent has finite rank $n,$ it can be written as an orthogonal sum of $n$ minimal tripotents (cf.\ \cite[Comments after Lemma 3.3 in page 200]{Ka97} and the previous Lemma \ref{l minimal tripotents in inner ideals} to deal with minimality).\label{finite rank of tripotents}
	
	It is well known that a C$^*$-algebra $A$ is reflexive as a Banach space if and only if it is finite dimensional. Actually, dimension and rank are intrinsically related in the setting of C$^*$-algebras. More concretely, every finite dimensional C$^*$-algebra must have finite rank by basic structure theory (cf.\ \cite[Theorem I.11.2]{TakBook78}). On the other hand, every infinite-dimensional C$^*$-algebra must have infinite rank (see, for example, \cite[Exercise 4.6.13]{KadRingVolI}). If $A$ is an infinite-dimensional von Neumann algebra, we can clearly find an infinite sequence of mutually orthogonal non-zero projections.
	Contrary to what is known in the case of C$^*$-algebras and projections, the rank and not the dimension of the Peirce-2 subspace associated with a tripotent $e$ in a JB$^*$-triple $E$ determines the number of mutually orthogonal elements in $E_2(e)$. We can find tripotents $e$ whose Peirce-2 subspace is infinite dimensional with finite rank. For example, every spin factor $C$ has rank 2 (cf.\ \cite{Ka97}). If $C$ is infinite-dimensional, it contains many unitary tripotents, and for each one of them the Peirce-2 subspace is the whole $C$, and thus it is infinite-dimensional, reflexive, and does not contain more than 2 mutually orthogonal tripotents. This particularity of JB$^*$-triples produces a subtle gap in the proofs of \cite[Theorems 3.5 and 3.6]{Santos2014} and \cite[Theorem 4.5]{MarOikh2004}, where it is assumed that for a tripotent $e$ in a JB$^*$-triple $E$, the hypothesis dim$(E_2(e))= \infty$ implies that its bidual contains as many mutually orthogonal tripotents as desired whose sum is $e$. We shall see in the proof of Theorem \ref{t ADPr for general wk polynomials} an argument to fix this small inconvenient.
	
	The next technical lemma is thought to clarify the application of the notion of rank in JB$^*$-triples.
	
	\begin{lemma}\label{l tripotent of infinite rank} Let $e$ be a tripotent in a JBW$^*$-triple $M$. Suppose that $M_2(e)$ is a JBW$^*$-triple having infinite rank. Then there exists an infinite sequence $(e_n)$ of mutually orthogonal non-zero tripotents in $E_2(e)$ with $\displaystyle e = \hbox{w$^*$-}\sum_{n=1} ^{\infty} e_n $.
	\end{lemma}
	
	\begin{proof} We begin by observing that $e$ cannot be minimal in $M$, otherwise, $M_2(e) = \mathbb{C} e$ which is impossible because $\mathbb{C} e$ has rank one.
		
		Since minimality in the set of tripotents of a JBW$^*$-triple is precisely being minimal with respect to the natural partial order (cf.\ \cite[Corollary 4.8]{EdRutt88} and \cite[Lemma 4.7]{Batt91}), we can write $e$ as the orthogonal sum of two non-zero tripotents $e_1$ and $e_2$. If $e_1$ and $e_2$ are minimal in $M_2(e)$, by Lemma \ref{l minimal tripotents in inner ideals}, $e_1$ and $e_2$ are minimal tripotents in $M$, and thus, by \cite[Lemma 2.7]{FerMarPe2004} or \cite[Lemma 3.6]{KalPe2021}, $M_2(e)$ would be reflexive and hence of finite rank, leading to a contradiction.
		
		We can therefore assume that $e_1$ is not minimal in $M_2(e)$ (equivalently, in $M$). Therefore, we can again write $e_1$ as the orthogonal sum of two non-zero tripotents $e_{11}$ and $e_{12}$, which produces a decomposition of $e$ as the orthogonal sum of three non-zero tripotents.
		
		We shall prove, by induction, that for each natural $n$ there exist non-zero mutually orthogonal tripotents $e_1, \ldots, e_n$ satisfying $e_j \leq e$ for all $j$ and at least one of them is not minimal in $M_2(e)$ (equivalently, in $M$). The cases $n=1$ and $n=2$ have been proved above. Suppose, by the induction hypothesis, the existence of non-zero mutually orthogonal tripotents $e_1, \ldots, e_n$ satisfying $e= e_1+ \cdots + e_{n}$ and at least one of them is not minimal in $M$. We can assume that $e_n$ is not minimal. Therefore, we can find two non-zero orthogonal tripotents $e_{n1}$ and $e_{n2}$ in $M$ such that $e_n = e_{n1} + e_{n2}$. Then $  e= e_1+ \cdots + e_{n-1}+ e_{n1} + e_{n2}.$ We shall prove that one of them is not minimal. Otherwise, $e$ writes as an orthogonal finite sum of minimal tripotents, and hence $e$ must have finite rank by the comments in page \pageref{finite rank of tripotents}, which contradicts our hypotheses. Finally the induction argument gives the desired conclusion.
	\end{proof}

	Pe\l{}czy{\'n}ski's Property $(V)$ is a powerful tool to study weak compactness. Let us recall that a series $\displaystyle\sum_{n\geq 1} x_n$ in a Banach space $X$ is called weakly unconditionally Cauchy (wuC) if $\displaystyle \sum_{n\geq 1} |\varphi(x_n)|<\infty$ for each $\varphi \in X^*$. As shown by Bessaga and Pe\l{}czy{\'n}ski, wuC series are essentially linked to the canonical basis of $c_0$ \cite[\S VI]{DiesBook84}. According to \cite{Pelc}, a  Banach space $X$ satisfies property $(V)$ if for any (bounded) non relatively weakly compact set $K\subset X^*$ there exists a wuC series $\displaystyle\sum_{n\geq 1} x_n$ in $X$ such that $\sup_{\varphi\in K} |\varphi (x_n)|$ does
	not converge to $0$. If $X$ is a Banach space satisfying property $(V)$, then a bounded linear operator $T$ from $X$ into any other Banach space $Y$ is either weakly compact or fixes an isomorphic copy of $c_0$ (cf.\ \cite{Pelc}). Consequently, each bounded linear operator from $X$ into a Banach space not containing $c_0$ is weakly compact. It is further known that the dual space of a Banach space having property $(V)$ is weakly sequentially complete (see \cite{Pelc}). It follows that the dual space of a Banach space satisfying property $(V)$ cannot contain a copy of $c_0$, and if $X$ and $Y$ satisfy property $(V)$, every bounded linear operator $T\colon X\longrightarrow Y^*$ is weakly compact.
	
	The list of Banach spaces satisfying Pe\l{}czy{\'n}ski's Property $(V)$ includes all C$^*$-algebras \cite[Corollary 6]{Pfitzner94}, and the strictly wider class of JB$^*$-triples \cite{ChuMe97}. Consequently, every bounded linear operator from a JB$^*$-triple $E$ into the dual space, $F^*,$ of any other JB$^*$-triple, $F,$ is weakly compact. The latter conclusion can be strengthened to the fact that every bounded linear operator $T\colon E\longrightarrow F^*$ actually factors through a Hilbert space (cf.\ \cite[Lemma 5]{ChuIoLou}).
	
	The procedure initiated by R.~Arens in \cite{Ar} for bilinear operators was later extended by R.~M.~Aron and P.~D.~Berner in \cite{ArBer78}, materializing in the most employed tool to produce norm-preserving extensions multilinear operators to the bidual spaces. Let $T\colon X_{1}\times X_{2}\times \cdots \times X_{N}\longrightarrow Y$ be a bounded multilinear map (where $X_{1},X_{2}, \ldots, X_{N},Y$ are Banach
	spaces). The Aron-Berner procedure guarantees the existence of $N!$ norm-preserving ``natural'' extensions of $T$ to multilinear maps from $X_{1}^{\ast \ast }\times X_{2}^{\ast \ast }\times \ldots \times X_{N}^{\ast \ast }$ to $Y^{\,\ast \ast }$, each one enjoying certain weak$^*$ continuity property. These Aron-Berner extensions are separately weak$^*$-to-weak$^*$ continuous if and only if they all coincide. A result by F.~Bombal and I.~Villanueva proves that, assuming that for each $i \neq j$, every bounded linear operator from $X_i$ into $X^*_j$ is weakly compact, every bounded multilinear operator $T\colon X_{1}\times X_{2}\times \cdots \times X_{N}\longrightarrow Y$ admits a unique norm-preserving separately weak$^*$-to-weak$^*$ continuous Aron-Berner extension $\widehat{T}\colon X_{1}^{\ast \ast }\times X_{2}^{\ast \ast }\times \ldots \times X_{N}^{\ast \ast } \longrightarrow Y^{\,\ast \ast }$ \cite[Theorem 1]{BomVi}.
	
	Another basic principle of functional analysis asserts that a bounded linear operator $T\colon X\longrightarrow Y$ (where $X$ and $Y$ are Banach spaces) is weakly compact if and only if $T^{**} (X^{**})\subseteq Y$ (cf.\ \cite[Theorem 5.5]{ConwayBook}). The setting of multilinear operators differs from the linear case. One implication remains true. Namely, in order to guarantee the uniqueness of the Aron-Berner extension, let us assume that for each $i \neq j$ every bounded linear operator from $X_i$ into $X^*_j$ is weakly compact. Then it is known that each weakly compact multilinear mapping $T\colon X_{1}\times X_{2}\times \cdots \times X_{N}\longrightarrow Y$ satisfies that its Aron-Berner extension is $Y$-valued  (cf.\ \cite[Corollary 2]{BomVi}). However, for $N>1$ there exist non-weakly compact $N$-linear operators whose Aron-Berner extension remains valued in the same codomain space (see, for example, \cite[page 385]{Pelc63} or \cite[Lemma 7]{GonGu95} and \cite[Theorem 5]{Vi}).
	
	In the case of multilinear maps from the Cartesian product of a collection of JB$^*$-triples $E_1, \ldots, E_N$ into a complex Banach space $X,$ the precise characterization of those satisfying that their Aron-Berner extension remains valued in the same codomain space was established by A.~M.~Peralta, I.~Villanueva, J.~D.~M.~Wright, and K.~Ylinen in \cite{PeViWriYl2010multi}. We recall that a multilinear operator $T\colon E_1\times \cdots \times E_N\longrightarrow X$ is \emph{quasi completely continuous} if whenever we choose strong$^*$ Cauchy sequences $(x_n^i)_n\subset E_i$ ($1\leq i \leq N$), it follows that the sequence $(T(x_n^1,\ldots, x_n^N))_n$ is norm convergent, equivalently, given sequences $(x^i_n)\subset E_i$
	which are strong$^*$ convergent to $x^i\in E_i$ ($1\leq i \leq k$), we have
	\begin{equation*}
		\lim_n \|T(x^1_n,\ldots, x^N_n)-T(x^1,\ldots,x^k)\|=0.
	\end{equation*}
	
	\begin{theorem}[\textrm{\cite[Theorem 3.9]{PeViWriYl2010multi}}]\label{t characterization PVWY multi} Let $E_1,\ldots,E_N$ be JB$^*$-triples, $X$ a complex Banach space, and $T \colon E_1\times	\cdots \times E_N \longrightarrow X$ a bounded multilinear operator. The following assertions are equivalent:
		\begin{enumerate}[$(1)$]
			\item $T$ is quasi completely continuous.
			\item The unique Aron-Berner extension of $T$ is $X$-valued.
		\end{enumerate}
	\end{theorem}
	
	A former version of the above theorem for C$^*$-algebras was obtained by J.~D.~M. Wright and K.~Ylinen (see \cite[Corollary 3.6]{WriYli}), while the case of abelian C$^*$-algebras was treated by I.~Villanueva in \cite{Vi} (see also \cite[Corollary 3.7]{WriYli}).

	Let $E$ be a JB$^*$-triple. Each scalar polynomial $p\colon E\longrightarrow \mathbb{C}$ writes as a finite sum $p = p_0+p_1+\cdots +p_N$, where each $p_k$ is a scalar $k$-homogenous polynomial. Clearly, Theorem \ref{t characterization PVWY multi} above implies that the symmetric multilinear form $\widehat{p}_k$ is quasi completely continuous, and hence $p_k\colon E\longrightarrow \mathbb{C}$ is sequentially strong$^*$(-to-norm) continuous for all $0\leq k\leq N$. The conclusion in the next corollary just follows by gluing the different $k$-homogeneous summands expressing $p$.
	
	\begin{corollary}\label{c scalar polynomials on JBstar triples as sequentilly strong* continuous} Every scalar polynomial on a JB$^*$-triple is sequentially strong$^*$ continuous.
	\end{corollary}
	
	\section{Main results}\label{section:proofs}
	
	We can now state and prove that every JB$^*$-triple satisfying the Daugavet property also satisfies the stronger polynomial Daugavet property.

	\begin{theorem}\label{t DPr for general wk polynomials} Let $E$ be a JB$^*$-triple satisfying the Daugavet property. Then $E$ satisfies the polynomial Daugavet property, that is, every weakly compact polynomial $P\colon E \longrightarrow E$ satisfies the Daugavet equation \eqref{eq DP}.
	\end{theorem}
	
	\begin{proof} Let $P \colon E \longrightarrow E$ be a  weakly compact polynomial. We can assume, via \cite[Theorem 1.1]{ChoGarMaestreMartin2007} (see also Proposition \ref{prop_charpolDP}), that $P$ is of the form $P(x) = p(x) a$, where $a\in E$, $p\colon E\longrightarrow \mathbb{C}$ is a scalar polynomial with $\|a\| = 1= \|q\|$. For each $\varepsilon>0$, we can find $x_0$ in the closed unit ball of $E$ and $\omega\in \mathbb{C}$ with $|\omega|=1$ such that $\re\omega p(x_0) > 1-\varepsilon.$ Let $\widehat{p} \colon E^{**}\longrightarrow \mathbb{C}$ denote the Aron-Berner extension of the polynomial $p$, obtained by extending each homogeneous summand of $p$.
		
		As in the argument in the proof of \cite[Theorem 3.5]{Santos2014} (or \cite{BeMar2005,MarOikh2004}), the key information is in the triple spectrum of $a$. Let $Sp(a) \subseteq [0,1]$ denote the triple spectrum of $a$ in $E$, where $1\in Sp(a)$.
		
		Our next goal will consist in proving that for each $\varepsilon>0$ there exists a sequence $(e_n)_n$ of mutually orthogonal tripotents in $E^{**}$ such that \begin{equation}\label{eq existence of the sequence of mo tripotents}\hbox{$P_2(e_n) (a) = \lambda_n e_n$ with $\lambda_n >1-\varepsilon,$ for all $n\in \mathbb{N}$.}
		\end{equation}
		
		If $1$ is isolated in $Sp(a)$, the element $e = \chi_{\{1\}}$ --i.e., the characteristic function corresponding to the set $\{1\}$-- is a continuous function in $C_0(Sp(a))\equiv E_a$ and lies in the JB$^*$-subtriple $E_a$ generated by $a$. Furthermore, $a = e + a_0$ with $a_0\perp e$. We shall next show that $E_2(e)$ has infinite rank. Otherwise it must be reflexive (cf.\ \cite[Proposition 4.5]{BuChu92} and \cite[Theorem 6]{ChuIo90} or \cite{BeLoRo03, BeLoPeRo04}). It is known that every reflexive JB$^*$-triple is generated by its minimal tripotents (cf.\ \cite[Proposition 4.5 and Remark 4.6]{BuChu92}). However, every minimal tripotent in $E_{2}(e)$ is a minimal tripotent in $E$ (cf.\ Lemma \ref{l minimal tripotents in inner ideals}) but the latter contains no minimal tripotents by hypothesis (cf.\ Theorem \ref{theorem:DPr-ADP-CstarJBstar}). This concludes the proof of the fact that $E_2(e)$ has infinite rank. By observing that $E_2(e)$ is a weak$^*$-dense JB$^*$-subtriple of $E^{**}_2(e),$ we conclude that the latter has finite rank (equivalently, is reflexive) if and only if  $E_2(e)$ has finite rank. Therefore  $E_2^{**}(e)$ must have infinite rank.

		Now, by applying Lemma \ref{l tripotent of infinite rank} to the JBW$^*$-triple $E^{**}$ and the tripotent $e$, we deduce the existence of an infinite sequence $(e_n)_n$ of mutually orthogonal non-zero tripotents in $E_2^{**}(e)$ with $\displaystyle e = \hbox{w$^*$-}\sum_{n=1} ^{\infty} e_n $. Let us recall that $a = e + a_0$ with $a_0\perp e$, therefore, by orthogonality, we get $$P_2(e_n) (a) = P_2(e_n) (e) = e_n = \lambda_n e_n, \hbox{ with }  \lambda_n = 1 > 1-\varepsilon.$$ This proves the statement in \eqref{eq existence of the sequence of mo tripotents} in the case in which $1$ is isolated in $Sp(a)$.
		
		If $1$ is non-isolated in $Sp(a)$, by \cite[Proposition 3.6]{ArKaup} we can find an infinite collection of mutually orthogonal non-zero tripotents $(e_n)_n$ in $E^{**}$ satisfying \eqref{eq existence of the sequence of mo tripotents}.
		
		We have already proved \eqref{eq existence of the sequence of mo tripotents}. Now, having in mind that $e_n \perp P_0 (e_n) (x_0)$ and \eqref{eq M-orhtogonality}, we define the following sequence in the unit sphere of $E^{**}$ by the rule $$x_n := \overline{\omega} e_n + P_0 (e_n) (x_0).$$ We note that $$x_0- x_n  = P_1(e_n) (x_0) + P_2(e_n) (x_0-  \overline{\omega} e_n) = P_1(e_n) (x_0) + P_2(e_n) (x_0) -  \overline{\omega} e_n.$$ The sequence $(e_n)_n$ converges to zero in the strong$^*$-topology of $E^{**}$ (see comments in page \pageref{eq orthogonal sequences are strong* null}), and the same occurs to $( \overline{\omega} e_n)_n$. We deduce from the joint strong$^*$ continuity of the triple product on bounded sets of $E^{**}$ that the sequences $(L(e_n,e_n)(x_0))_n = (P_2(e_n)(x_0) + \frac12 P_1(e_n) (x_0))_n,$ $(Q(e_n)(x_0))_n$ and $(P_2(e_n)(x_0))_n = (Q(e_n)^2(x_0))_n$ converge to zero in the strong$^*$-topology of $E^{**}$. In particular, $(P_2(e_n)(x_0))_n,$ $(P_1(e_n) (x_0))_n \longrightarrow 0$ in the strong$^*$-topology of $E^{**}$. This guarantees that $(x_n- x_0) \longrightarrow 0$ and $(x_n) \longrightarrow x_0$ in the strong$^*$-topology of $E^{**}$.
		
		By Corollary \ref{c scalar polynomials on JBstar triples as sequentilly strong* continuous} the polynomials $p\colon E\longrightarrow \mathbb{C}$ and $\hat{p} \colon E^{**}\longrightarrow \mathbb{C}$ are sequentially strong$^*$ continuous. Therefore $\hat{p} (x_n) \longrightarrow \hat{p}(x_{0})= p(x_0)$. It is now time to pick $n_0$ in $\mathbb{N}$ such that $\re\omega \hat{p}(x_{n_0}) > 1- \varepsilon.$ Our choice of $n_0$ implies the following facts: $P_2(e_{n_0}) (x_{n_0}) = \overline{\omega} e_{n_0}$ and $$\begin{aligned} \left\| x_{n_0} + \hat{p} (x_{n_0}) a \right\| & \geq \left\| P_2(e_{n_0}) \left(x_{n_0} + \hat{p} (x_{n_0}) a \right)\right\|
			= \left\| \overline{\omega} e_{n_0} + \hat{p} (x_{n_0})  \lambda_{n_0} e_{n_0}\right\| \\
			&= \left| \overline{\omega} + \hat{p} (x_{n_0})  \lambda_{n_0}  \right| = \left| 1 + {\omega} \hat{p} (x_{n_0})  \lambda_{n_0}  \right| \\
			&\geq \re \left( 1 + {\omega} \hat{p} (x_{n_0})  \lambda_{n_0}  \right) > 2 - 2 \varepsilon,
		\end{aligned}$$ witnessing that $\left\|\id_{E} + P \right\| = \left\|\id_{E} - {p} a\right\| = \left\|\id_{E^{**}} - \widehat{p} a\right\| > 2 - 2 \varepsilon$, and the arbitrariness of $\varepsilon$ gives the desired conclusion.
	\end{proof}	
	
	Let us now deal with the alternative Daugavet Property.
	
	\begin{theorem}\label{t ADPr for general wk polynomials} Let $E$ be a JB$^*$-triple satisfying the alternative Daugavet property. Then $E$ satisfies the alternative polynomial Daugavet property, that is, every weakly compact polynomial $P\colon E \longrightarrow E$ satisfies the alternative Daugavet equation \eqref{eq ADP}.
	\end{theorem}

	\begin{proof} Since $E$ satisfies the ADP every minimal tripotent in $E$ is diagonalizing (see Theorem \ref{theorem:DPr-ADP-CstarJBstar}). Suppose $P\colon E\longrightarrow E$ is a weakly compact polynomial. We can assume, via \cite[Corollary 1.2]{ChoGarMaestreMartin2007} (see Proposition \ref{prop_charpol-ADP}), that $P(x) =  p(x) a$, where $p\colon E\longrightarrow \mathbb{C}$ is a polynomial and $a\in X$. There is no loss of generality in assuming that $\|p\| = 1= \|a\|$.
		
		Let us refine a bit our previous arguments in the proof of Theorem~\ref{t DPr for general wk polynomials}. If $1$ is isolated in the triple spectrum of $a$, the characteristic function $\chi_{_{\{1\}}}$ is a tripotent $e$ in $C_0(Sp(a))\equiv E_a$, where the latter denotes the JB$^*$-subtriple of $E$ generated by $a$.
		
		The dichotomy is the following: $E_2(e)$ has finite or infinite rank. In the first case $E_2(e)$ must be reflexive (cf.\ \cite[Proposition 4.5]{BuChu92} and \cite[Theorem 6]{ChuIo90} or \cite{BeLoRo03, BeLoPeRo04}), and the whole $E$ is generated by its minimal tripotents (\cite[Proposition 4.5 and Remark~4.6]{BuChu92}). Therefore, $E_2(e)$ contains a minimal tripotent $e_1$ which is also minimal in $E$ by Lemma \ref{l minimal tripotents in inner ideals}, and diagonalizing by hypotheses. Therefore $E$ decomposes in the form $E= \mathbb{C} e_1 \oplus_{\infty} E_0(e_1)$ ($E_2(e_1) = \mathbb{C} e_1$). We note that $e = \gamma e_1 + P_0(e_1) (e)$ for a suitable unitary $\gamma\in \mathbb{C}$.
		
		Since $\|p\|=1$, for each $\varepsilon>0$, there exists $\alpha\in \mathbb{C}$ with $|\alpha|\leq 1$ and $x_0$ in the closed unit ball of $E_0(e)$ such that $|p(\alpha e_1 + x_0)|> 1-\varepsilon$. Having in mind that $p$ is a scalar polynomial, the mapping $\zeta \longmapsto  p(\zeta \gamma e_1 + x_0)$ is a polynomial on $\zeta\in \mathbb{C}$, and thus an entire function. By the maximum modulus principle, there exists a unitary scalar $\lambda_1$ such that \begin{equation}\label{eq px1 big} |p(\lambda_1 \gamma e_1 + x_0)| = \max\bigl\{|p(\beta e_1 + x_0)| \colon \beta \in \mathbb{C}, |\beta|\leq 1\bigr\} \geq |p(\alpha e_1 + x_0)| > 1-\varepsilon.
		\end{equation} Let us observe that $a = e+ P_0 (e) (a)= \gamma e_1  + P_0(e_1) (a)$ because $\gamma e_1\leq e$. Setting $x_1 := \lambda_1 \gamma e_1 + x_0,$ we define a norm-one element in $E$ satisfying the following:
		$$\begin{aligned} \|(\id_{E} +\omega P)(x_1)\|&= \|x_1 +\omega p(x_1) a \| \geq \| P_2(e_1) (x_1 +\omega p(x_1) a  )\| \\ &= \| \lambda_1 \gamma e_1 + \omega p(x_1) \gamma e_1  \| =  | \lambda_1 + \omega p(x_1) | = | 1 + \omega \overline{\lambda_1} p(x_1) |,
		\end{aligned}$$ for every unitary $\omega\in \mathbb{C}$. It is straightforward to check, from the previous inequality, that $\max\limits_{|\omega|=1} \|(\id_{E} +\omega P)(x_1)\| = 1 + | p(x_1) |> 2- \varepsilon$ (see \eqref{eq px1 big}), witnessing that $\max\limits_{|\omega|=1} \|\id_{E} +\omega P\| = 2$ as desired.
		
		The above argument shows that every polynomial of the form $P (x) = p(x) a$ with $p\in \mathcal{P} (E),$ $a\in E,$ and $\|p\| = \|a\|=1$, satisfies the alternative Daugavet equation \eqref{eq ADP} whenever $1$ is isolated in the triple spectrum of $a$ and there exists a minimal tripotent $e_1$ in $E$ with $e_1 \in  E_2(e)$ (in particular, when $E_2(e)$ has finite rank).
		
		In the remaining cases one of the next statements holds:\begin{enumerate}[$(a)$] \item  $1$ is isolated in the triple spectrum of $a$ and for $e= \chi_{\{1\}}\in E$ and the Peirce-2 subspace $E_2(e)$ has infinite rank.
			\item $1$ is non-isolated in the triple spectrum of $a.$
		\end{enumerate}
		In both cases, we can literally repeat the arguments in the proof of Theorem \ref{t DPr for general wk polynomials} to show that $P$ satisfies the Daugavet equation and thus the alternative Daugavet equation.
	\end{proof}
	
	Obviously, the conclusions in Theorems \ref{t DPr for general wk polynomials} and \ref{t ADPr for general wk polynomials} actually characterize those JB$^*$-triples satisfying the Daugavet property and the alternative Daugavet property, respectively. We have indeed rediscovered the conclusions in Theorem \ref{theorem:DPr-ADP-CstarJBstar} because we have proved that $E$ has the DPr (respectively, the ADP) provided that $E$ contains no minimal tripotents (respectively, every minimal tripotent in $E$ is diagonalizing).
	
	We conclude this paper with an interesting open problem.
	
	\begin{problem}
		Let $E$ be a JB$^*$-triple satisfying the Daugavet property {\rm(}respectively, the alternative Daugavet property{\rm)}.
		In view of Theorem \ref{t characterization PVWY multi} {\rm(}\cite[Theorem 3.9]{PeViWriYl2010multi}{\rm)}, it seems natural to ask whether every polynomial $P\colon E \longrightarrow E$ whose Aron-Berner extension remains $E$-valued satisfies the \eqref{eq DP} {\rm(}respectively, the \eqref{eq ADP}{\rm)}.
	\end{problem}

	\smallskip\smallskip
	
	
	\textbf{Acknowledgements}

D.~Cabezas supported by Junta de Andaluc\'ia I+D+i grants P20\_00255 and FQM-375.
M.~Mart\'in supported by Project PGC2018-093794-B-I00 funded by MCIU/AEI/10.13039/501100011033/FEDER, UE, by Junta de Andaluc\'ia I+D+i grants P20\_00255, A-FQM-484-UGR18, and FQM-185, and by ``Maria de Maeztu'' Excellence Unit IMAG CEX2020-001105-M/AEI/10.13039/501100011033/FEDER. A.M.~Peralta supported by MCIN/AEI/10.13039/501100011033/FEDER ``Una manera de hacer Europa'' project no.\ PGC2018-093332-B-I00, Junta de Andaluc\'{\i}a grants FQM375 and  A-FQM-242-UGR18, and by the IMAG--Mar{\'i}a de Maeztu grant CEX 2020-001105-M/AEI/10.13039/501100011033.
	

\end{document}